\documentclass[conference]{IEEEtran}
\IEEEoverridecommandlockouts
\usepackage{cite}
\usepackage{amsmath,amssymb,amsfonts}
\usepackage{attrib}
\usepackage{xcolor}
\usepackage{array}
\usepackage{tikz}
\usetikzlibrary{positioning, arrows.meta}
\usepackage{algorithmic}
\usepackage{graphicx}
\usepackage{textcomp}
\usepackage{xcolor}
\def\BibTeX{{\rm B\kern-.05em{\sc i\kern-.025em b}\kern-.08em
    T\kern-.1667em\lower.7ex\hbox{E}\kern-.125emX}}

    \usepackage{fancyhdr}
\fancypagestyle{arxivnotice}{%
  \fancyhf{} % Clear all headers and footers
  \fancyfoot[C]{\parbox{\textwidth}{\footnotesize \copyright~2026 IEEE. Personal use of this material is permitted. Permission from IEEE must be obtained for all other uses, in any current or future media, including reprinting/republishing this material for advertising or promotional purposes, creating new collective works, for resale or redistribution to servers or lists, or reuse of any copyrighted component of this work in other works.}}
}

\begin{document}

\title{Design and Implementation of Interactive Videos in Undergraduate Transnational Education

}

\author{Dimitrios Chiotis\\
\IEEEauthorblockA{\textit{School of Physical and} \\ \textit{Chemical Sciences} \\
\textit{Queen Mary University of London}\\
London, United Kingdom \\
d.chiotis@qmul.ac.uk}
\and
 Ebru Burcu \\
\IEEEauthorblockA{\textit{School of Physical and} \\ \textit{Chemical Sciences} \\
\textit{Queen Mary University of London}\\
London, United Kingdom \\
e.burcu@qmul.ac.uk}}
% \and
% \and
% \and
% \and

\maketitle
\thispagestyle{arxivnotice}

\begin{abstract}
The growing development of online educational video resources over the past two decades, has significantly reshaped university students' self-regulated study approaches and influenced diverse ways of learning. Within an environment of abundantly and publicly available online resources, educators cannot always control which materials students choose to supplement their study. Consequently, considerable deviations are likely to be present between the taught content by the educator and an arbitrary online educational video the student watches, in terms of terminology, notation and methods, especially in mathematical sciences. These deviations can be even more considerable when Transnational Education (TNE) students incorporate native language resources to comprehend the material. This paper presents the design, creation and implementation of interactive video resources produced in alignment with established pedagogical frameworks to support an undergraduate TNE mathematics module. 
\end{abstract}


\section{Introduction}
Online videos have become a prominent component of university students’ learning practices, supporting activities ranging from entertainment to independent study and revision. Their use is particularly evident in mathematically intensive courses, where students frequently engage with step-by-step demonstrations and visualisations of abstract concepts through video resources \cite{KAY2012619}. Within the mathematical sciences, such resources can be especially valuable, as sequential argumentation and symbolic manipulation often require paced exposure that students can pause, revisit and review as needed \cite{inglis2012proof}. It is only natural that video-based learning is nowadays considered an appropriate complement to in-person lectures.  

In \cite{KAY2012619} it was further established that students utilise video resources to support learning of advanced concepts. Video materials have been proved to support mathematics fluency and conceptual understanding via the incorporation of worked examples, the revision of key steps in problem-solving and the exposition to multiple representations of mathematical ideas. Moreover, Wirth et al. \cite{wirth2024comparing} support that video explanations can positively affect student engagement with maths reasoning in proof-based contexts.  In the context of tablet-based handwritten recorded videos, Yoon and Sneddon \cite{yoon2011} report that students highly valued the ability to access material and review problem-solving steps at their own pace.  

Although publicly available learning video resources possess certain advantages, at the same time they pose certain academic challenges. According to Mayer and Florella's coherence principle of multimedia learning \cite{mayer2014coherence}, extraneous material, such as deviations from notation, terminology and problem-solving methods used in-course, can obscure learning by overloading working memory, especially to students who are still developing foundational understanding.  

These deviations can become even more apparent in TNE settings, where material in different languages might be selected by the students. 
As the viewer is not an expert in the studied field and encounters complex mathematical ideas from different viewpoints, inconsistencies of this kind can be carried forward and establish incorrect foundations in the academic development. In TNE contexts, these challenges may be further expanded, as students might use materials from different educational systems and languages. 

Instead of restricting students’ usage of online video resources, this paper suggests a more academic approach: educators can adapt by creating and providing structured videos which align with their taught content.

In this paper we describe the central workflow of creating interactive video resources to support an undergraduate mathematics course which forms part of a TNE programme. The videos combine tablet-based problem solving, lecturer video capturing and embedding of interactive elements to enhance self-reflection and self-assessment. Also, we outline the pedagogical framework underpinning the design of the videos and the editing considerations before publishing them to the module’s virtual learning environment (VLE).

\section{Pedagogical frameworks}
The development of video resources, which are presented in this paper, was pedagogically grounded in established principles in multimedia learning and mathematics education. The Cognitive Theory of Multimedia Learning provides a pedagogical foundation for designing teaching content which combines visual and verbal components. Principles such as signalling, personalisation and embodiment have proved to enhance student engagement while minimising additional cognitive load, as supported by \cite{mayer2009multimedia}. In addition, research findings on worked examples \cite{sweller2011cognitive} highlight the importance of structured guidance in steps, especially when supporting the development of problem-solving skills to new learners. 

Within mathematics education, the scaffolding approach is considered important in encouraging learners to engage with complex topics. This usually involves guided support, questioning and incremental progress to enhance comprehension of unfamiliar material \cite{bakker2015scaffolding}. In parallel, the Community of Inquiry framework emphasises the impact of social and teaching presence in shaping relevant learning experiences in online contexts \cite{garrison2000critical}. Together, these perspectives highlight the need for structured, guided and engaging learning resources that support both conceptual understanding and student participation.

Finally, the design of the video resources also supports independent learning by encouraging elements of self-regulated learning. In this context, learners are required to actively monitor their understanding, engage with partially completed solutions and revisit key stages of the material as needed. Such practices align with established models of self-regulated learning, which emphasise active control over the learning process, including planning, monitoring and self-evaluation, as well as the role of motivation in sustaining engagement \cite{zimmerman2002becoming, pintrich1999role}.

\section{Motivation}
The work presented in this paper was developed within a mathematics module of a TNE programme. 

The development of these videos was largely motivated by preliminary findings indicating our students’ strong inclination to supplement their study routines with online video resources. Specifically, 80\% of respondents reported that they regularly watch online videos to supplement their learning, with 72\%
expressing a preference of video duration between 5 and 30 minutes. Also, 60\% of respondents indicated that they use such videos to understand difficult or newly introduced concepts and review lecture content. However, the main concern raised from the respondents was the difficulty to find reliable and academically aligned online materials.  

Consequently, our work focuses on addressing these needs by developing video resources aligned with the academic quality standards established in the lectures. The videos intend to support students in comprehending newly taught material through a scaffolding approach, while maintaining consistency in terminology, notation and methodology with the in-person lectures. 

\section{Design of interactive videos}
We propose a structured design framework for interactive mathematics videos that integrates pedagogical principles with practical implementation.

The design was informed by the pedagogical principles outlined in Section II. The approach combined two formats with complementary roles: lecturer-presented videos, used to introduce the topic, and tablet-based videos, where the studied method was illustrated step-by-step. These two formats were designed to provide conceptual explanation and guided problem-solving. 

 The purpose of the videos was to supplement and complement class-taught material and guide the students through the Simplex Method. Our design was based on three core principles: (i) structuring of mathematical content through guided step-by-step problem solving, (ii) inclusion of instructor presence to support engagement, (iii) incorporation of interactive elements to promote active learning. These design principles reflect established approaches in multimedia learning, scaffolding and self-regulated learning, linking conceptual explanation with structured guidance and active learner engagement.

 Lecturer-presented recordings included clear articulation of learning objectives, concept recall and mathematical motivation. Tablet-based recordings focused on elicitation and guided questioning, encouraging students to actively engage with the material and test their understanding at different stages. Also, tablet-based recordings adopted a scaffolding approach, with gradually minimised support to guide learners towards independent problem solving. Moreover, common mistakes were addressed aimed at strengthening conceptual understanding. 

\subsection{Lecturer-presented videos} 
These recordings were designed to introduce and frame each topic, aligning with the Community of Inquiry by establishing teaching presence from the start. The primary aim of these recordings was to define the learning objectives and provide the mathematical motivation underpinning the Simplex Method.

These videos contextualised new material by linking  it with previously taught topics, highlighting both its relevance within the module and its broader mathematical significance. The recordings incorporated clear articulation of learning objectives alongside signalling and concept recall, activating prior knowledge and supporting continuity in learning. This approach aimed at providing students with sufficient intuition and conceptual grounding before engaging with the technical details of the method.  

\begin{figure}[hbtp]
    \centering
    \includegraphics[width=\linewidth]{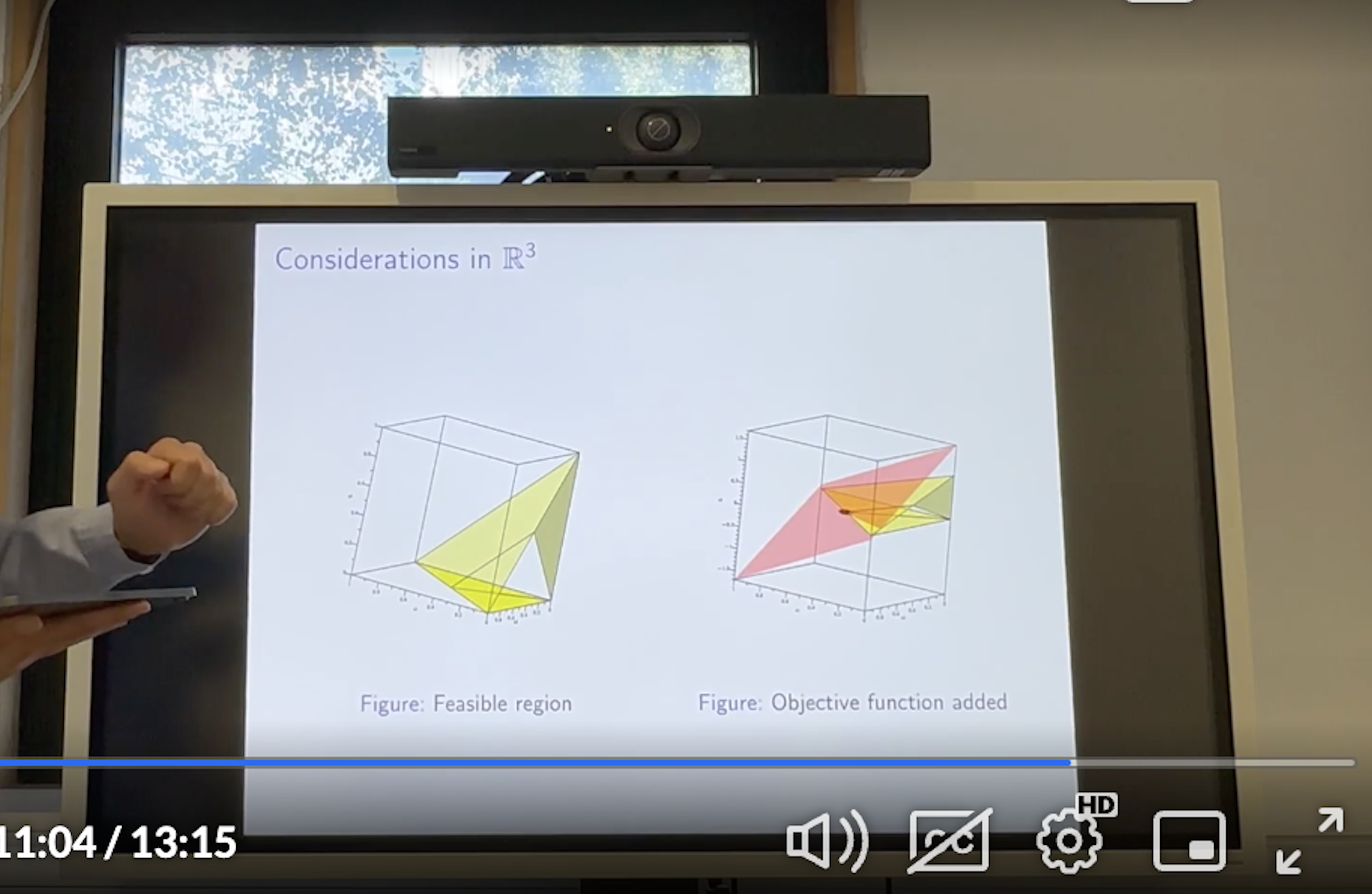}
    \caption{Lecturer-presented video reinforcing teaching presence.}
    \label{fig:lecpres}
\end{figure}

The inclusion of a visible lecturer presence  contributed to a more personalised learning experience. This was achieved through a combination of screen-based presentation and live whiteboard writing, with the former being shown in Fig. \ref{fig:lecpres}, allowing both structured exposition and incremental development of mathematical reasoning. These elements intended to enhance engagement, support clarity of explanation and foster a sense of connection with the material, particularly in independent study contexts. 

\subsection{Tablet-based problem solving} 
Tablet-based recordings were the central component of the video resources and focused on step-by-step development of mathematical derivations and solutions, aligning with worked example approaches which enhance the development of problem-solving skills in new learners. This format was specifically selected to illustrate the procedural aspects of the Simplex Method using digital ink technologies, allowing students to follow the reasoning in a sequential manner. Each stage of the derivation and solution was intentionally explicit in order to model mathematical thinking in a comprehensive and accessible manner.  

\begin{figure}[!b] \centering\includegraphics[width=\columnwidth]{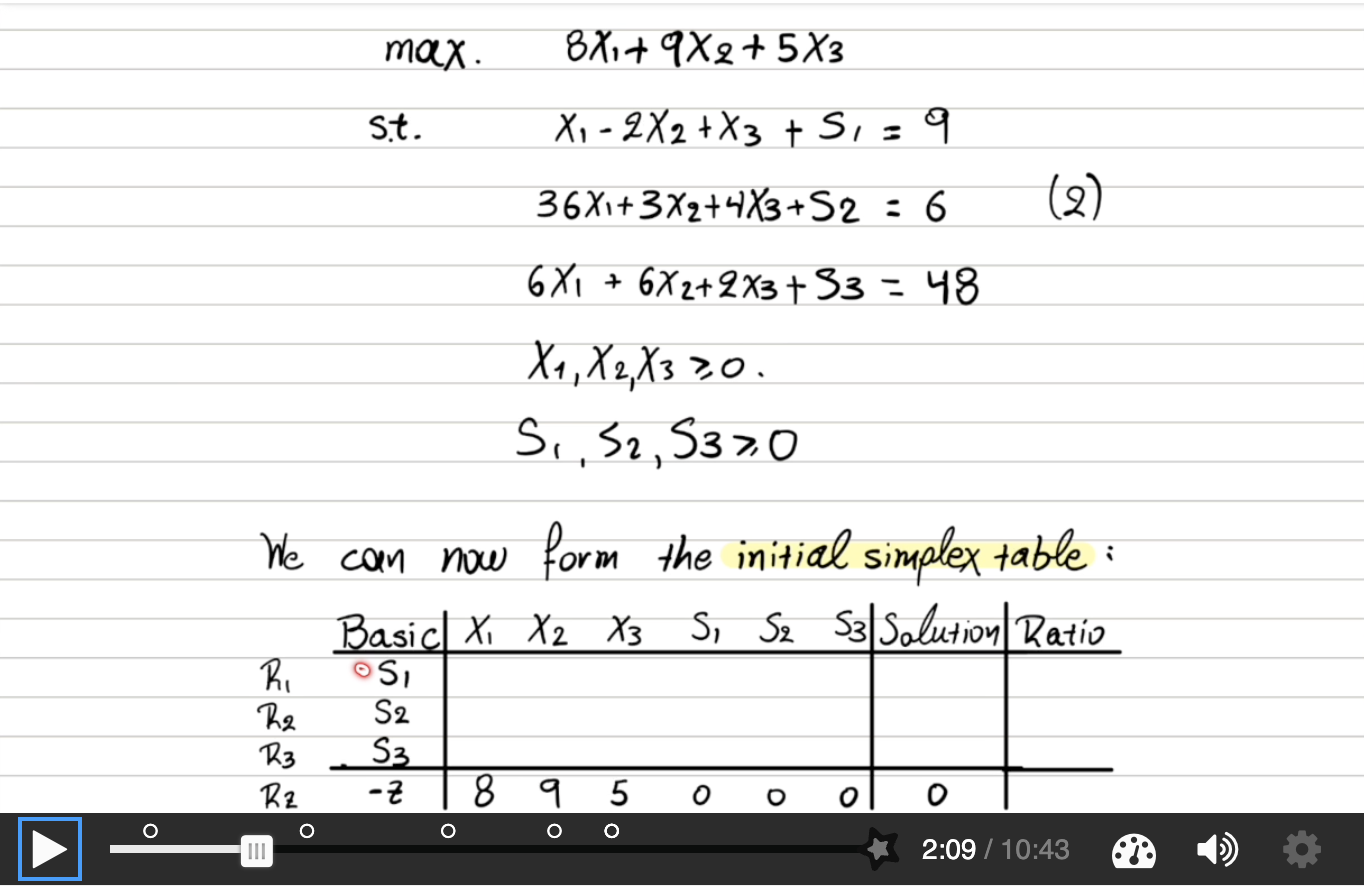}
     \caption{Partially completed step within a tablet-based video, prompting students to actively engage by completing the required table before proceeding.}
     \label{fig:placeholder1}
 \end{figure}

To support progressive learning, two tablet-based videos for each variation of the Simplex Method were created: standard cases, unbounded linear programs and multiple optimal solutions. The first video in each pair provided a full step-by-step guided demonstration of the method, and the second required more of a student input (via interactive elements, see below) with certain steps intentionally omitted, as in Fig. \ref{fig:placeholder1}, to encourage viewers to test their understanding by completing them independently. This progression aligns with scaffolding learning design, encouraging the development of incremental problem-solving abilities.

The aforementioned paired structure was designed to encourage active participation and support a gradual transition from passive demonstration to more active and independent problem-solving. With the second video containing partially completed answers, students were required to apply the presented methods, enriching both procedural fluency and conceptual understanding of the Simplex Method. This approach supports aspects of self-regulated learning, requiring students to monitor their understanding and actively engage with the omitted steps.

\subsection{Interactive elements}
\begin{figure}[!t]
    \centering
    \includegraphics[width=\linewidth]{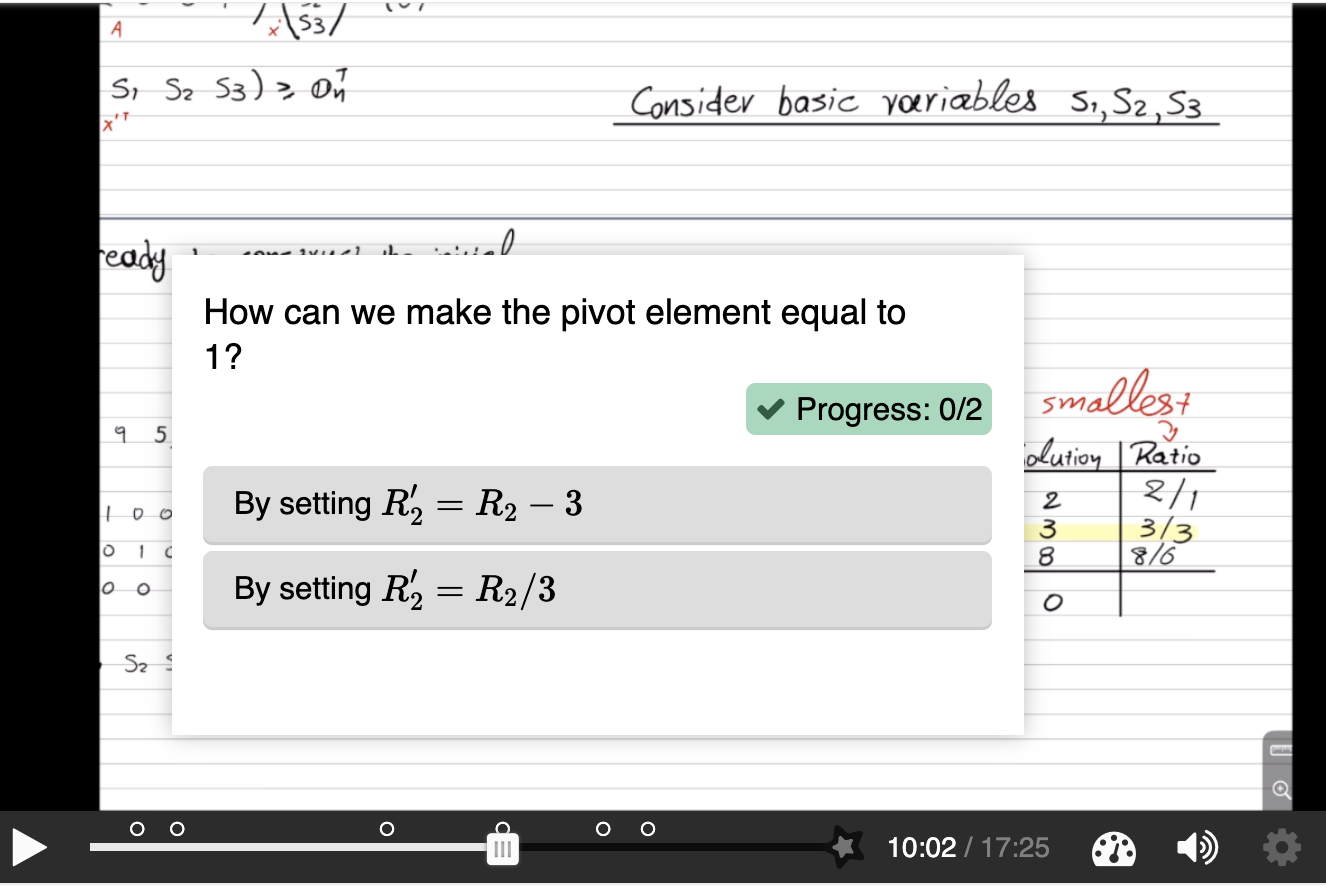}
    \caption{Interactive checkpoint embedded within a tablet-video, allowing students to engage with the material by answering a question before proceeding.}
    \label{fig:placeholder}
\end{figure}

These were incorporated within both video resources in each pair. They were embedded as pop-up checkpoints at key stages, especially following the introduction of new concepts and techniques in problem-solving, as shown in Fig. \ref{fig:placeholder}. 

The interactive elements were visible along the video play bar, allowing students to identify and revisit them, and also navigate directly to specific checkpoints. 
The pop-up prompts required students to reflect on the material and apply the new ideas before resuming the video. Multiple choice questions were employed and the ability to use LaTeX formatting in VLE helped ensure consistency and alignment with the notation used in lectures and textbooks.

These elements were designed to interrupt passive viewing and encourage active engagement, supporting processes associated with self-regulated learning, such as monitoring understanding and revisiting key material.

\subsection{Video editing}
Video editing was performed locally using computer native software, where both lecturer-presented and tablet-based recordings were refined to improve clarity, pacing and overall quality. Interactive elements were  later added using Lumi software, enabling the  integration of checkpoints and questions within the videos.

\subsection{Publishing to the VLE}
The completed interactive video resources were exported as H5P files from Lumi and were deployed to the VLE via Kaltura, with interactive elements carried through to the VLE. This allowed students to engage with embedded checkpoints and questions directly within the module's VLE.
\subsection{Practical Considerations}
The development and embedding of interactive video resources in the VLE requires consideration of both pedagogical and technical factors. From a production perspective, the level of preparation needed may vary depending on the experience and confidence of the lecturer. In particular, the use of written prompts or structured outlines can support lecturers in delivering coherent explanations when recording both camera-based presentations and tablet-based problem-solving sessions. Such preparation can contribute to a more natural delivery while maintaining consistency in the presentation of mathematical content.

In addition,  video editing and the integration of interactive elements require a basic level of technical proficiency. While the tools used in this work were accessible, the editing and deployment stages may be perceived as time-consuming or demanding for some users, particularly those with limited prior experience. Hence, adequate support and familiarisation with the required software may be necessary to ensure an efficient implementation.

\section{Conclusion }
The proposed design highlights the feasibility and potential of integrating structured video resources within TNE education. The combination of conceptual explanation, step-by-step problem solving and embedded interaction supports both active engagement and independent learning practices. 

In practice, the effectiveness of this approach depends on how clearly the material is structured and how well the interactive elements are incorporated. The design also requires a level of technical adaptation, which may present challenges for some educators.

This paper provides a practical example of how pedagogically informed video resources can be developed and embedded within a module. Future work will focus on evaluating the impact on student performance and engagement.

\bibliography{references}
\end{document}